\documentclass{agtart_a}
\pdfoutput=1

\usepackage{pinlabel}

\title[Concordance invariants of doubled knots]
{Ozsv\'{a}th--Szab\'{o} and Rasmussen invariants\\of doubled knots}

\author{Charles Livingston}
\givenname{Charles}
\surname{Livingston}
\address{Department of Mathematics\\
  Indiana University\\
  Bloomington, IN 47405\\USA}
\email{livingst@indiana.edu}

\author{Swatee Naik}
\givenname{Swatee}
\surname{Naik}
\address{Department of Mathematics and Statistics\\
University of Nevada\\
Reno, NV 89557\\USA}
\email{naik@unr.edu}
\thanks{Research supported by the NSF}

\volumenumber{6}
\issuenumber{}
\publicationyear{2006}
\papernumber{24}
\startpage{651}
\endpage{657}

\doi{}
\MR{}
\Zbl{}


\received{26 February 2006}
\revised{}
\accepted{2 March 2006}
\published{18 May 2006}
\publishedonline{18 May 2006}
\proposed{}
\seconded{}
\corresponding{}
\editor{}
\version{}

\arxivreference{math.GT/0505361}

\AtBeginDocument{\let\bar\wbar\let\tilde\wtilde\let\hat\what}

\newcommand{\fig}[2]{\includegraphics[scale=#1]{\figdir/#2}}

\makeatletter
\def\cnewtheorem#1[#2]#3{\newtheorem{#1}{#3}
\expandafter\let\csname c@#1\endcsname\c@theorem}

\newtheorem{theorem}{Theorem}
\cnewtheorem{lemma}[theorem]{Lemma}
\cnewtheorem{corollary}[theorem]{Corollary}

\theoremstyle{definition}
\cnewtheorem{definition}[theorem]{Definition}

\makeatother  

\numberwithin{equation}{section}

\begin{document}

\begin{asciiabstract}
Let \nu be any integer-valued additive knot invariant that bounds the
smooth 4-genus of a knot K, |\nu(K)| \le g_4(K), and determines the
4-ball genus of positive torus knots, \nu(T_{p,q}) = (p-1)(q-1)/2.
Either of the knot concordance invariants of Ozsvath-Szabo or
Rasmussen, suitably normalized, have these properties.  Let
D_{\pm}(K,t) denote the positive or negative t-twisted double of K.
We prove that if \nu(D_{+}(K,t)) = \pm 1, then \nu(D_{-}(K,t)) = 0.
It is also shown that \nu(D_{+}(K,t))= 1 for all t \le TB(K) and
\nu(D_{+}(K, t)) = 0 for all t \ge -TB(-K), where TB(K) denotes the
Thurston-Bennequin number.

A realization result is also presented: for any 2g \times 2g Seifert
matrix A and integer a, |a| \le g, there is a knot with Seifert form
A and \nu(K) = a.
\end{asciiabstract}

\begin{htmlabstract}
<p class="noindent">Let &nu; be any integer-valued additive
knot invariant that bounds the smooth 4&ndash;genus of a knot K,
|&nu;(K)|&le;g<sub>4</sub>(K), and determines the 4&ndash;ball genus of
positive torus knots, &nu;(T<sub>p,q</sub>)=(p-1)(q-1)/2.  Either of
the knot concordance invariants of Ozsv&aacute;th&ndash;Szab&oacute;
or Rasmussen, suitably normalized, have these properties.
Let D<sub>&plusmn;</sub>(K,t) denote the positive or
negative t&ndash;twisted double of K. We prove that if
&nu;(D<sub>+</sub>(K,t))=&plusmn;1, then &nu;(D<sub>-</sub>(K,t))=0. It
is also shown that &nu;(D<sub>+</sub>(K,t))=1 for all t&le;TB(K) and
&nu;(D<sub>+</sub>(K,t))=0 for all t&ge;-TB(-K), where TB(K) denotes
the Thurston&ndash;Bennequin number.</p>

<p class="noindent">A realization result is also presented: for any
2g&times;2g Seifert matrix A and integer a, |a|&le;g, there is a knot
with Seifert form A and &nu;(K)=a.</p>
\end{htmlabstract}

\begin{abstract}
Let $\nu$ be any integer-valued additive   knot invariant that bounds the
smooth 4--genus of a knot $K$,  $|\nu(K)| \le  g_4(K) $, and determines
the 4--ball genus of positive torus knots, $\nu(T_{p,q}) = (p-1)(q-1)/2$.
Either of the knot concordance invariants of Ozsv\'{a}th-Szab\'{o}   or
Rasmussen, suitably normalized, have these properties.  Let $D_{\pm}(K,t)$
denote the positive or negative  $t$--twisted double of $K$. We prove
that if $\nu(D_{+}(K,t) )= \pm 1$, then  $\nu(D_{-}(K,t)) = 0$. It is
also shown that $\nu(D_{+}(K,t) )= 1$ for all $t \le \mbox{TB}(K)$
and  $\nu(D_{+}(K, t) )= 0$ for all $t \ge -\mbox{TB}(-K)$, where
$\mbox{TB}(K)$ denotes the Thurston-Bennequin number.

A realization result is also presented: for any $2g \times 2g$ Seifert
matrix $A$ and integer $a$, $|a| \le g$, there is a   knot   with Seifert
form $A$ and $\nu(K) = a$.
\end{abstract}

\maketitle

\section{Introduction}  

Two recently discovered smooth knot concordance invariants, the
Ozsv\'{a}th-Szab\'{o} invariant $\tau$,~\cite{os}, and the Rasmussen
invariant $s$,~\cite{ra}, have opened up powerful new perspectives on
the study of concordance.  For instance, each is sufficient to prove
the Milnor conjecture determining the smooth 4--ball genus of torus
knots, and each demonstrates the existence of non-slice Alexander
polynomial one knots.  Unfortunately, the computation of $\tau$ is not
algorithmic and $s$, though algorithmic, is difficult to compute for
classes of knots. The only classes for which general results are known
are alternating knots and torus knots
(Ozsv\'{a}th-Szab\'{o}~\cite{os}, Rasmussen~\cite{ra}), and
quasipositive knots (Livingston~\cite{liv}, Shumakovich~\cite{sh}).
Related results are included in Eftekhary, Hedden--Ording,
Owens--Strle and Plamenevskaya \cite{ef,ho,ow,pl}.  Of special
relevance here are the recent results for doubles of $(2,n)$--torus
knots by Matt Hedden and Philip Ording~\cite{ho}, discussed further in
an addendum to this paper.

The three observations of this note grew out of efforts to compute
$\tau$ and $s$.  In~\cite{liv}, methods for computing $\tau$ were
developed and used to prove that some untwisted doubles of knots have
non-vanishing $\tau$.  The results of~\cite{liv} depended on only
three properties of $\tau$: (1) $\tau(K \# J) = \tau(K) + \tau(J)$,
(2) $|\tau(K)| \le g_4(K)$, and (3) $\tau(T_{p,q}) =(p-1)(q-1)/2$ for
all torus knots $T_{p,q}$ with $p, q >0$.  The invariant $s$ shares
these properties of $\tau$ when suitably normalized; more precisely,
$-s/2$ has these three properties. Thus, for this paper we refer to an
arbitrary knot invariant having these three properties as $\nu$.

Let $D_{\pm}(K,t)$ denote the positive or negative  $t$--twisted double of $K$.  Since $D_{\pm}(K,t)$ bounds a surface of genus 1, $|\nu( D_{\pm}(K,t) | \le 1$.  Our results are the following.

\begin{theorem}\label{thmd1}   If $\nu(D_{+}(K,t) )= \pm 1$, then  $\nu(D_{-}(K,t)) = 0$.  
\end{theorem}

 Associated to a knot $K$ there is an integer-valued invariant called
 the Thurston-Bennequin number, TB($K$), first defined in
 Bennequin~\cite{b}. This invariant was initially defined in terms of
 Legendrian structures, but there is the following combinatorial
 definition.  Every knot has some diagram so that in a neighborhood of
 each crossing the projection consists of two segments, one with slope
 1 and the other of slope $-1$.  Furthermore it can be arranged that
 the segment with slope 1 passes under the other segment.  For such a
 diagram, the Thurston-Bennequin number is the writhe of the diagram
 minus the number of right cusps (maximum points with respect to
 projection onto the $x$--axis). The Thurston-Bennequin number of the
 knot is the maximum value of this difference, taken over all diagrams
 satisfying the crossing criteria.

 For further details about the Thurston-Bennequin number, see Ng~\cite{n}.   Note that $\mbox{TB}(-K) \ne -\mbox{TB}(K)$; in fact, for all $K$,   $\mbox{TB}(-K) +  \mbox{TB}(K) \le -1$.

\begin{theorem}\label{thmd2} For each knot $K$ there is an integer $t_K$ such that    
$\nu(D_+(K,t)) = 1$ for $t \le t_K$ and  $\nu(D_+(K,t)) = 0$ for all $t > t_K$.  The value of $t_K$ satisfies $ \mbox{TB}(K) \le t_K < -\mbox{TB}(-K)$
\end{theorem}

 Similar results hold with the roles of $D_+$ and $D_-$ reversed. 
 
Techniques related to those we use here to study doubled knots also
quickly yield the following.

\begin{theorem}\label{thmr1} For any $2g \times 2g$ 
Seifert matrix $A$ and integer $a$, $|a| \le g$, there is a    knot  $K$ bounding a genus $g$ Seifert surface with Seifert form $A$ and $\nu(K) = a$.
\end{theorem}

{\bf Acknowledgments}\qua We thank Matt Heddon for sharing the results
of his work with Philip Ording with us.  We also thank Ciprian
Manolescu and Brendan Owens for identifying an error in sign
conventions in an early version of this paper.

 \section{Band modifications and 4--genus}
\begin{figure}[h]
\begin{center}
\labellist
\pinlabel {\tiny $t$} at 184 603
\pinlabel {\tiny $K$} at 312 618
\endlabellist
\fig{.35}{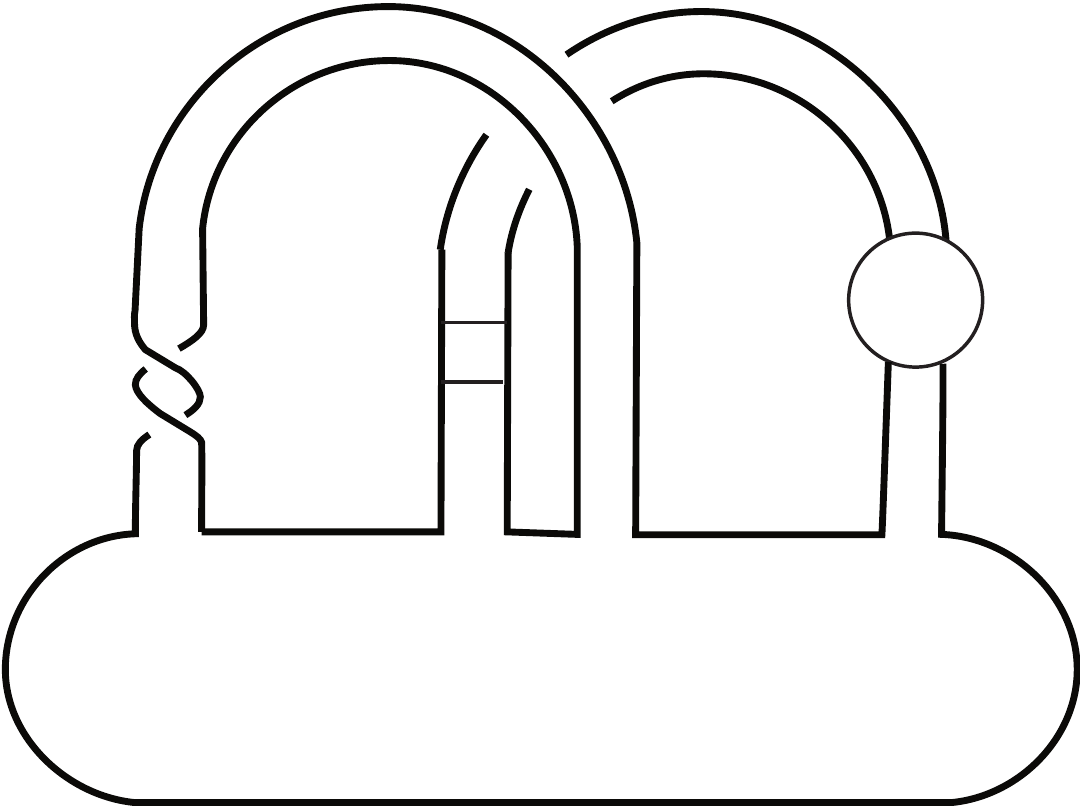} \hskip.07in 
\labellist
\pinlabel {\tiny $t$} at 377 370
\pinlabel {\tiny $K$} at 504 379
\endlabellist
\fig{.35}{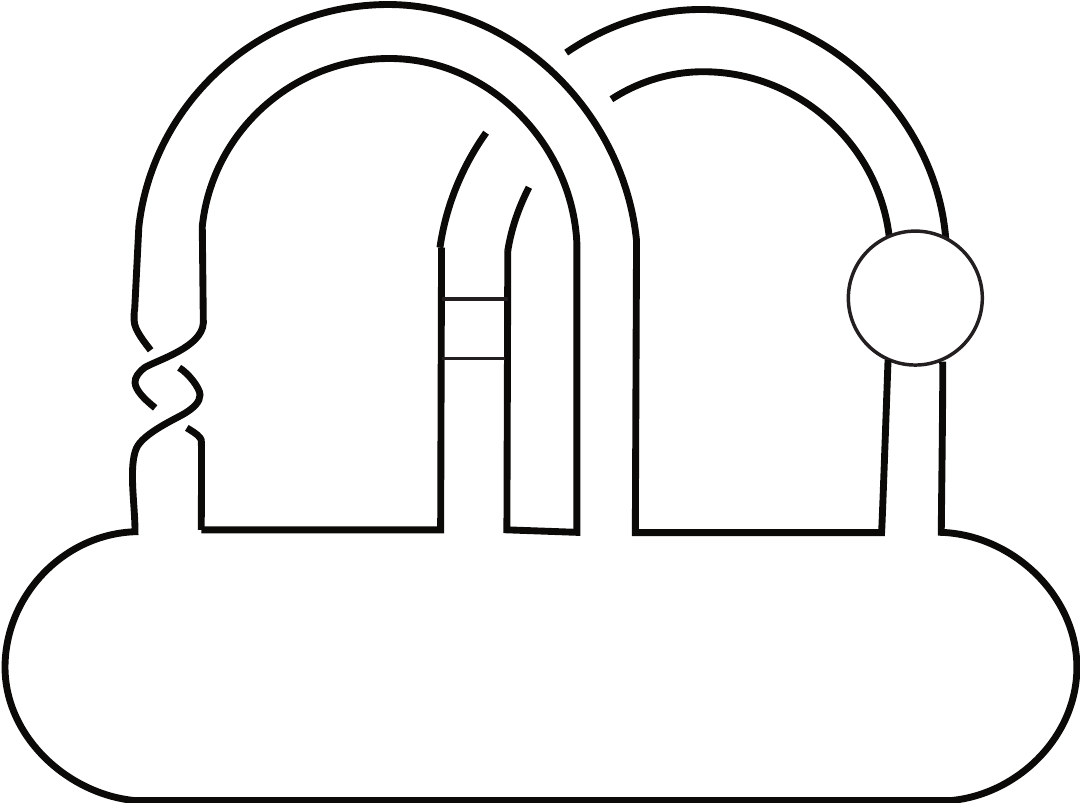} \hskip.07in 
\labellist
\pinlabel {\tiny -$t$} at 156 132
\pinlabel {\tiny --$K$} at 284 149
\endlabellist
\fig{.35}{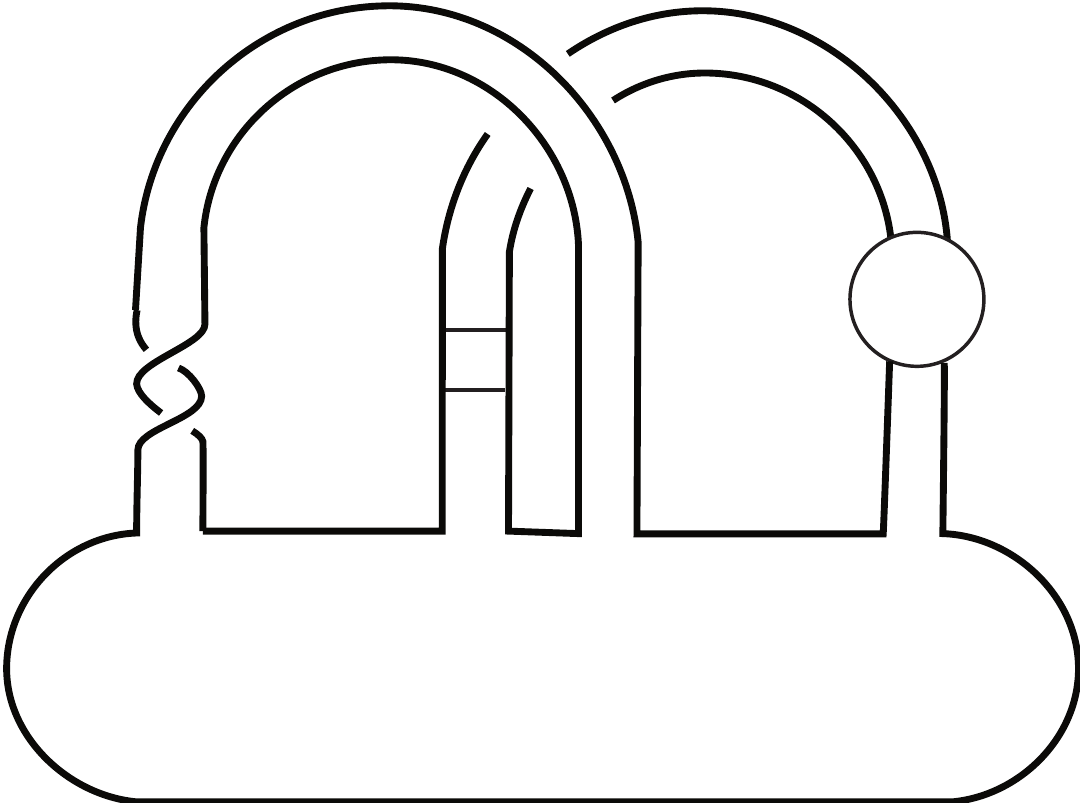}
\end{center}
\caption{$D_+(K,t),\hskip.1in D_-(K,t), \hskip.0in \  \mbox{and}\   D_-(-K,-t) = -D_+(K,t) $} \label{figure1}
\end{figure}

A Seifert surface for a knot, viewed as a disk with bands attached, can be modified by removing one of the bands and reattaching it in the same place,  perhaps twisted and knotted  in a different way.  We call this operation   {\it band modification}. For example, the first diagram in \fullref{figure1} illustrates the knot $D_+(K,t)$: the band on the right has framing $t$.  A band modification is performed on the left band to yield $D_-(K,t)$.  A band modification is then performed on the right band to yield $D_-(-K,-t)$,  and we note that $ D_-(-K,-t)= -D_{+}(K,t)$.  

\begin{theorem} \label{genthm} If a knot $K'$ is constructed from a knot $K$ by performing a single band modification, then $K \# - K'$ has smooth 4--genus 0 or 1.
\end{theorem}
\begin{proof} For any knot $J$, removing a single band in the Seifert surface yields a 2--component  link $L$, and there is a genus 0 cobordism from $J$ to $L$.  In our setting removing the corresponding bands from   $K$ and $K'$ yields the same link $L$.  Gluing the two cobordisms together along $L$ yields a genus 1 cobordism from $K$ to $K'$.  The standard construction of removing the neighborhood of an arc from $K$ to $K'$ on that cobordism yields a punctured surface of genus 1, bounded by $K \# - K'$, embedded in $B^4$.
\end{proof}

We assume now  and for the rest of the paper that $\nu$ is an integer-valued additive knot invariant satisfying $|\nu(K)| \le g_4(K)$ for all $K$ and $\nu(T_{p,q}) = (p-1)(q-1)/2$ for all torus knots $T_{p,q}$ with $p, q > 0$.

\begin{corollary}\label{corm} If $K'$ is constructed from $K$ by a band modification, then $|\nu(K) - \nu(K')| \le 1$.
\end{corollary}
\begin{proof} Since $K' \# -K'$ is slice, $|\nu(K') + \nu(-K')| = |\nu(K' \# -K')| = 0$ and so  $ \nu(-K')=  -\nu(K')$.  By \fullref{genthm} we then have $|\nu(K)  - \nu(K')| = |\nu(K \# -K')| \le 1$.  
\end{proof} 

\section{Knot invariants and  doubles}\label{doublessection}

\noindent{\bf \fullref{thmd1}}\qua {\sl If $\nu(D_+(K,t)) = \pm 1$ then $\nu(D_-(K,t)) = 0$.}  
 
\begin{proof}  We have that $\nu(-D_+(K, t)) = \mp 1$.  From the  example  illustrated in \fullref{figure1}, we see that $D_-(K,t)$ can be  constructed from   $D_+(K,t)$ and also from $-D_+(K,t)$ by a single band modification.  Thus, by \fullref{corm},  $\nu(D_-(K,t))$ differs from both $1$ and $-1$ by at most 1, so it must be 0.
\end{proof}

\noindent{\bf \fullref{thmd2}}\qua {\sl For each knot $K$ there is an integer $t_K$ such that    
$\nu(D_+(K,t)) = 1$ for $t \le t_K$ and  $\nu(D_+(K,t)) = 0$ for $t > t_K$.  The value of $t_K$ satisfies $ \mbox{TB}(K) \le t_K < -\mbox{TB}(-K)$.}

\begin{proof}  The construction here is much the same as one originated by Rudolph~\cite{ru} as formulated  in~\cite{liv}; we only summarize it here.  

For any $t \le \mbox{TB}(K)$, $K$ can be isotoped to be on a minimal genus Seifert surface for some positive torus knot, $T_{p,q}$, with self-framing $t$.    It  follows that a genus one Seifert surface for $D_+(K,t)$ embeds on a torus knot Seifert surface. By results of~\cite{liv}  one then has $\nu(D_+(K,t) )=   1$.

Similarly, for   $t \le \mbox{TB}(-K)$, $\nu(D_+(-K,t)) = 1$. By \fullref{thmd1}, it follows that $\nu(D_-(-K,t)) = 0$.  Taking the mirror image, $\nu(D_+(K,-t)) = 0$.  Equivalently, $\nu(D_+(K,t)) = 0$ if $t \ge -\mbox{TB}(-K)$.

As proved in~\cite{liv}, the three conditions satisfied by $\nu$ imply that changing a positive crossing to a negative crossing in a knot diagram cannot increase the value of $\nu$ (and can decrease it by at most 1).  The knot $D_+(K,t)$  results from $D_+(K,t+1)$ when a negative crossing in the knotted band is changed into a positive crossing; notice that the crossings in a positively twisted band are negative crossings, since the strands are oriented in opposite directions.  Thus, $\nu(D_+(K,t))$ is a non-increasing function of $t$.  The integer $t_K$ is defined to be the largest integer for which $\nu(D_+(K,t)) =1$.  Its bounds follow from the results of the previous two paragraphs.  
\end{proof}

{\bf Example}\qua  For the  right-handed trefoil we have TB($ T_{2,3}) = 1$ and TB($- T_{2,3}) = -6$ (see~\cite{n}).  Thus,  $\nu(D_{+}(T_{2,3},t)) = 1$ for $t \le 1$ and   $\nu(D_{+}(T_{2,3},t)) = 0$ for $t \ge 6$.  (See  the addendum for remarks on recent work of Matt Hedden and Philip Ording~\cite{ho} concerning the values of $\tau$ and $s$ for  twisted doubles of $(2,n)$--torus knots.) 

\noindent{\bf Remark}\qua The asymptotic limiting behavior of $\nu$ for
doubled knots holds for more general families of companions.  Suppose
that $J$ is a winding number 0 knot in the solid torus that meets a
meridinal disk transversely in exactly two points.  Given any knot
$K$, one can form the $t$--twisted satellite, $J(K,t)$, using $J$ as
the companion.  As $t$ is increased, $J(K,t)$ is changed by adding
negative crossings.  Thus, the value of $\nu(J(K,t))$ is
non-increasing by a result of~\cite{liv}.  On the other hand, the
genus of $J(K,t)$ is bounded above by the genus of a surface bounded
by $J$ in the solid torus, so $\nu(J(K,t))$ is bounded below by the
negative of that genus.  It follows that the value of $\nu(J(K,t))$
has some finite limit as $t$ increases.

In the case that $J$ intersects the meridinal disk in more than 2
points, this argument does not work.  Further properties of $\tau$ as
it relates to surfaces in $CP^2$, proved in~\cite{os}, can be applied
to recover the asymptotic behavior in this case, but no similar
argument is known for $s$.

\section{Realization result}
 We now restate and prove the realization result:
 
\medskip
{\bf \fullref{thmr1}}\qua
 {\sl
Given a $2g \times 2g$ Seifert matrix $A$,  
there are  
knots $ K_i,\ -g \le i \le g$, bounding Seifert surfaces of genus $g$ 
each with  $A$ as its associated Seifert matrix,
such that $\nu (K_i) \, =\, i$. 
}

\begin{proof}  It is essentially an observation of Rudolph, the Trefoil 
Insertion Lemma in~\cite{ru}, that by repeatedly adding trefoils to
the bands of the Seifert surface in such a way that the framings of
the bands, and thus the Seifert form, is unchanged, one eventually
arrives at a surface $S$ that is isotopic to an embedded surface on a
minimal genus Seifert surface for some positive torus knot.  According
to~\cite{liv}, for $K = \partial S $,
$$\nu (K)\, =\, g_4(K)\, =\, g_3(K)\, =\, g(S).$$

Apply this construction for the Seifert matrix $A$ to build a   knot $K_+$ that bounds a genus $g$ Seifert surface $S_+$  having Seifert form $A$ and $\nu(K_+) = g$.  Similarly, build a    knot $K_-$ with genus $g$ Seifert surface  $S_-$, Seifert form $-A$ and $\nu(K_-) = g$.  Then $K_+$ and $-K_-$  both bound genus $g$ Seifert surfaces ($S_+$ and $-S_-$, respectively) with Seifert form $A$, where   $\nu(K_+) = g$ and $ \nu(-K_-) = -g$.

  A series of $2g$ band modification converts $S_+$ into $-S_-$: just replace the bands for $S_+$ with those of $-S_-$ one at a time, maintaining the Seifert form at each step, as described in the following paragraph.  Since each of the modifications changes the value of $\nu$ by at most 1, and the $2g$ modifications decrease $\nu$ by $2g$, it follows that each step reduces $\nu$ by exactly one.  Thus the sequence of knots arising from the series of modifications yields the desired examples.
  
  To conclude the proof, we indicate the process of trading bands from those of the first surface to those of the second while maintaining the Seifert form. Given knots $K$ and $K'$ with the same Seifert form, we can assume both are built from  the same disk $D$ by adding bands, $\{b_1, \ldots b_{2g}\}$ and  $\{b'_1, \ldots b'_{2g}\}$.  We can further assume that the attaching maps for $b_i$ and $b'_i$ are the same, and that all the bands are disjoint. The band $b_1$ can be removed and replaced with $b'_1$.  However, the linking numbers associated to  $b'_1$ with the other bands can differ from that of $b_1$.  This is corrected by stretching $b'_1$ to add clasps between  it and the   $b_i, i \ge 2$.  Using an isotopy, $b'_1$ can be returned to its original position, at the expense of moving the $b_i, i \ge 2$.  This procedure can now be repeated, with $b_2$ and $b'_2$, and so on, until the desired result is achieved.

\end{proof}

{\bf Addendum (March, 2006)}\qua Since this paper
 was circulated, Hedden and Ording~\cite{ho} have proved that the
 Ozsv\'ath-Szab\'o invariants and the Rasmussen invariant are
 distinct.  Let $t_\tau(K)$ and $t_s(K)$ denote the greatest value of
 $t$ for which $\tau(D_+(K,t)) = 1$ or $s(D_+(K,t))= -2$, a value which
 exists by the results of this paper.  Hedden and Ording compute
 $t_\tau(T_{2,3}) = 1$ but $t_s(T_{2,3}) \ge 2$.  More generally, they
 show that $t_\tau(T_{2,2n+1}) = 2n -1$ and computations of Rasmussen
 (usng Bar-Natan's program~\cite{ba}) show that $t_s(T_{2,5}) \ge 5$,
 and $t_s(T_{2,7}) \ge 8$.

\bibliographystyle{gtart}
\bibliography{link}

\end{document}